\documentclass{elsart}

\usepackage{amssymb}
\journal{J. Math. Anal. Appl.}

\newtheorem{theorem}{Theorem}[section]

\theoremstyle{definition}

\newtheorem{example}[theorem]{Example}

\newtheorem{remark}[theorem]{Remark}

\newcommand*{\pd}[2]{\frac{\partial #1}{\partial #2}}
\newcommand*{\pdd}[2]{\frac{\partial ^2 #1}{\partial #2 ^2}}

\begin{document}
\begin{frontmatter}




\title{On the dependence of the reflection operator on boundary conditions for
biharmonic functions}
\author[Tanya]{Tatiana  Savina \corauthref{cor}}
\thanks[bby]{Research of this author was supported in part by OU
Research Challenge Program, award \# RC-09043.  } \ead{savin@ohio.edu}
\ead[url]{http://www.math.ohiou.edu/ tanya/}
\corauth[cor]{Corresponding author.}
\address[Tanya]{Department of Mathematics, Condensed Matter \& Surface Science Program, Nanoscale \& Quantum Phenomena Institute, Ohio University, Athens,
OH 45701, USA}

\begin{abstract}
{The biharmonic equation arises in areas of continuum mechanics
including linear elasticity theory and the Stokes flows, as well as
in a radar imaging problem. We discuss the reflection formulas for
the biharmonic functions $u(x,y)\in\mathbb{R}^2$ subject to
different boundary conditions on a real-analytic curve in the plane.
The obtained formulas, generalizing the celebrated Schwarz symmetry
principle for harmonic functions, have different structures. In
particular, in the special case of the boundary, $\Gamma _0
:=\{y=0\}$, reflections are point to point when the given on $\Gamma
_0$ conditions are $u=\partial _nu=0$, $u=\Delta u=0$ or $\partial
_nu=\partial _n\Delta u=0$, and point to a continuous set when
$u=\partial _n\Delta u=0$ or $\partial _nu=\Delta u=0$ on $\Gamma
_0$.}
\end{abstract}

\begin{date}
\date{}
\end{date}

\end{frontmatter}
\maketitle

\section{  Introduction  }

In this paper we study the dependence of the structure of the
reflection operator on different boundary conditions for the biharmonic  functions,
  where a
function $u(x, y)$ of class $C^{4}(U)$ is said to be biharmonic   if
it is a solution to the equation $\Delta^2u=0$ \cite{lipkin}. Here
$U$ is a domain in $\mathbb{R}^2,$ and $\Delta$ is the  Laplacian.


The obtained reflection operator is generally an integro-differential operator,
which reduces in the simplest case to the celebrated local (point to
point) Schwarz symmetry principle for harmonic functions.

 In the case of harmonic functions, there are  three basic types of
boundary conditions: the
Dirichlet, Neumann and Robin, and the Schwarz reflection principle
can be stated as follows.
\begin{theorem}
Let $\Gamma \subset \mathbb{R}^2$ be a non-singular real analytic
curve and $P'\in \Gamma$. Then, there exists a neighborhood $U$ of
$P'$ and an anti-conformal mapping $R:U\to U$ which is identity on
$\Gamma$, permutes the components $U_1, U_2$ of $U\setminus \Gamma$
and relative to which any harmonic function $u(x, y)$ defined near
$\Gamma $ and \\
-- vanishing on $\Gamma $ (the homogeneous Dirichlet condition) is
odd \cite{dima}, \cite{shapiro},
\begin{equation} \label{E:0.1}
         u(x_0 ,y_0 ) =-u(R(x_0 ,y_0 )),
\end{equation}
-- subject to the Neumann condition on $\Gamma$, $\pd{u}{n}=0$, is
even,
\begin{equation} \label{E:0.1b}
         u(x_0 ,y_0 ) =u(R(x_0 ,y_0 )),
\end{equation}
-- subject to the Robin  condition on $\Gamma $,
$\alpha (x,y) \pd{u(x,y)}{n} +\beta (x,y) u(x,y)=0,\;(x,y)\in
\Gamma, $
can be continued by the following integro-differential operator \cite{be-sa},
\begin{equation}\label{E:1.2}
u(x_0 ,y_0 )=u(R(x_0 ,y_0 ))
\end{equation}
$$+\frac{1}{2i}\int\limits_{\Gamma }^{R(x_0 ,y_0 )}F(x,y,x_{0},y_{0})\omega(u(x,y))-u(x,y)
\omega(F(x,y,x_{0},y_{0})),$$ for any point $(x_0, y_0)$
sufficiently close to $\Gamma$. Here  $\omega (\cdot )
=\pd{}{y}dx-\pd{}{x}dy$, and the integral is independent on the path
joining an arbitrary point on $\Gamma$ with the point $R(x_0 ,y_0
)$. The function $F$ is defined by the coefficients $\alpha (x,y)$
and $\beta (x,y)$, and the curve $\Gamma$ \cite{be-sa}. The mapping
$R$ mentioned above is given by
\begin{equation} \label{E:0.2}
R(x_0 ,y_0 )=R(z_0 )=\overline{S(z_0 )},
\end{equation}
where $S(z)$ is the Schwarz function \cite{davis}.
\end{theorem}
Note that if the point $(x_0 ,y_0)\in U_1$, then the ``reflected''
point $R(x_0 ,y_0)\in U_2$, and the mapping $R$ depends only on the
curve $\Gamma$ and is defined only near $\Gamma$ but may have
conjugate-analytic continuation to a larger domain.

The Schwarz reflection principle has been studied by several
researchers (see \cite{as}, \cite{be-sa}, \cite{bramble}--\cite{G}, \cite{john}--\cite{lopez},
\cite{nystedt}--\cite{shapiro} and references therein). In
particular, when $\Gamma$ is a line, H.~Poritsky \cite{poritsky}
proved that a biharmonic function $u(x,y)$,  a solution to the
biharmonic equation $\Delta _{x,y}^2u=0$, defined for $y\ge 0$ and
subject to conditions
\begin{equation}\label{Dirichlet1}
u(x,0)=\pd{u}{y}(x,0)=0
\end{equation}
can be continued across the x-axis using the formula
\begin{equation} \label{E:0.02}
 u(x_0 ,y_0 ) =-u(R(x_0 ,y_0 )) -2y_0 \pd{u}{y}(R(x_0 ,y_0 ))
-y_0 ^2 \Delta _{x,y} u(R(x_0 ,y_0 )),
\end{equation}
where $R(x_0 ,y_0 )=(x_0 ,-y_0 )$ and
$\Delta _{x,y}=\pdd{}{x}+\pdd{}{y}$.
He also applied this formula to problems of planar elasticity, studying bending
of plates, where $u$ is a deflection of a thin plate, clamped along $y=0$.
An analogous formula has been obtained by R.J.~Duffin \cite{duffin} for
three-dimensional case. Duffin also considered spherical boundaries
and applied his result to study viscous flows, among other things.
J.~Bramble \cite{bramble} considered continuation of biharmonic
functions across the circular arc with the clamped boundary  conditions
\begin{equation}\label{Dirichlet2}
u=\pd{u}{r}=0 \quad \mbox{for} \quad x^2 +y^2 =a ^2,
\end{equation}
He has shown that
$u$ can be continued using the formula
\begin{eqnarray} \label{E:0.04}\nonumber
 u(x_0 ,y_0 ) =  -u(R(x_0 ,y_0 ))& \\
  - \frac{r^2 _0 -a ^2}{r^2 _0}
\bigl ( r_0 & \pd{u}{r}u(R(x_0 ,y_0 ))+\frac{1}{4}(r^2 _0 -a ^2 )
\Delta _{x,y} u(R(x_0 ,y_0 )) \bigr  ) ,\nonumber
\end{eqnarray}
where $r_0 =\sqrt{x^2 _0 +y^2 _0}$ and  $a$ is the radius of the
circle. J.~Bramble also applied it to the elastic medium problems in
two- and three-dimensional space. Papers by F.~John \cite{john} and
L.~Nystedt \cite{nystedt} are devoted to further studies of
reflection of solutions of linear partial differential equations
with various linear conditions on a hyperplane.
R.~Farwig in \cite{farwig} considered reflection principle for
biharmonic functions subject to different boundary conditions on a
hyperplane. Following to H.A.~Lorenz  and R.J.~Duffin \cite{duffin},
R.~Farwig \cite{farwig} applied his results to the Stokes system.

The purpose of this paper is to derive and to  compare  reflection formulas for
biharmonic functions across real analytic curves in $\mathbb{R}^2$
subject to  different boundary conditions occurring in
a variational approach to $\Delta ^2$. 
The motivation of this study is as follows, most of the results for biharmonic functions mentioned above were obtained for the Dirichlet boundary conditions, see (\ref{Dirichlet1}) and (\ref{Dirichlet2}). However, unlike the Laplace equation, which is typically complemented with either the Dirichlet,  Neumann or their linear combination, the Robin conditions,
 the biharmonic equation has much larger set of possible boundary conditions
as it is discussed in details in \cite{GGS}.
For this study we have chosen the conditions
(i) $u=\partial _n u =0$, (ii) $u=\Delta u=0$, (iii) $u=\partial _n \Delta u =0$
(iv) $\partial _nu=\Delta u=0$, (v) $\partial _n u=\partial _n\Delta u=0$ and (vi) $\Delta u=\partial _n\Delta
 u=0$,  
appearing  when one is considering a bilinear form associated with 
the biharmonic differential operator for $u,v\in C^4(\bar\Omega )$ \cite{mikhlin}, \cite{necas}
$$
a(u,v):=<\Delta ^2u,v>=\int\limits _{\Omega} v\Delta ^2u \,d\omega .
$$
Employing  Green's formula to $a(u,u)$, we obtain
\begin{equation}\label{firstGreen}
a(u,u)=
\int\limits _{\partial\Omega} (u\, \partial _n\Delta u-\partial _n u \, \Delta u) \, ds +\int\limits _{\Omega} (\Delta u)^2\, d\omega.
\end{equation}
The integral along the boundary $\partial\Omega$ in (\ref{firstGreen}) disappears under conditions (i), (ii), (v) and (vi), and therefore $a(u,u)\ge 0$. This immediately imply that the operator is strictly positive for (i) and (ii), and the corresponding boundary value problems in $\bar\Omega$ are well-posed. In the case  (v),  an extra condition, $\int\limits _{\Omega} u\,d\omega=0$, is usually imposed.  This yields the coercivity of the 
bilinear form in suitable function spaces.

Conditions (i), (ii) and (v) are commonly used in 
  applied problems for the biharmonic equation as well as for other equations with the biharmonic operator in the principal elliptic part. To give some examples, we remark that
condition (i), the Dirichlet condition,  corresponds to the clamped plate model   \cite{bramble}, \cite{duffin}, \cite{poritsky}, while  condition (ii), the Navier boundary condition,  corresponds to the hinged plate, when the contribution of curvature is neglected \cite{gazzola1}, \cite{GGS}. In the later case   the plate is ideally hinged along all of its edges so that it is free to rotate and does not experience any torque or bending moments about the edges.    Both the Dirichlet and Navier conditions are also used to model electrostatic actuation \cite{lin} alone with the governing non-linear nonlocal elliptic equation with the bi-Lanlacian in the principal part.

Condition (v) for biharmonic equation was considered in \cite{temam}, moreover, this condition is often used for the famous Cahn-Hilliard equation, $u_t=-\Delta (\Delta u+u-u^3)$, which is a semi-linear parabolic equation describing (among other processes) spinodal decomposition \cite{amy}, \cite{temam}. Pattern formation resulting this phase transition (spinodal decomposition) has been observed in alloys as well as polymer solutions and glasses. In this case the physical meaning of the second condition in (v), $\partial _n\Delta u=0$, is a no-flux condition (none of the mixture can pass through the walls of the container). The first condition in (v), $\partial _n u=0$, is the most natural way to ensure that the total free energy of the mixture decreases in time (this is a requirement from thermodynamics - the variational condition).

We have to comment on  condition (vi), for which we will not be
able to derive a new reflection formula. Indeed, consider biharmonic
function in $U$ subject to conditions (vi) on $\Gamma$. Denote
$v(x,y):=\Delta u(x,y)$, then function $v$ is a harmonic function in
$U$ subject to conditions $v=\partial _n v=0$ on $\Gamma$. Thus,
$v\equiv 0$ in $U$, and therefore $u(x,y)$ is a harmonic function in
$U$. Since any harmonic function is a solution to the problem in question, it means
that the space of solution does not have finite dimension.
Note also that as it was shown in \cite{GGS}, condition (vi) does not satisfy the Complementing Condition, introduced by S.~Agmon, A.~Douglis, and L.~Nirenberg \cite{agmon}. This condition is necessary for obtaining estimates up to the boundary for solutions of boundary value problems to the elliptic equations and, therefore,  is crucial for  the existence and uniqueness results \cite{agmon}. 

In the case of conditions (iii) and (iv) only one of the term in the boundary integral (\ref{firstGreen}) disappears, so our interest to this conditions is mostly theoretical, however,
 one of the problems arising in radar imaging may be reformulated as a boundary value problem for the biharmonic functions with  condition $\partial _n u=u -\beta\, \Delta u =0$ \cite{golub}, which coincide with (i) if $\beta =0$ and with (iv) if $\beta =-\infty$.
We remark, that other linear combinations, the Robin-type conditions, involving  (i)-(vi) are used 
in applications as well, but the corresponding reflections will be discussed elsewhere.


 The reflection formula  for the case of condition (i)
was obtained in \cite{as}. The only  known in the literature result
for the conditions (ii) and (v),   is in \cite{farwig}, where the
boundary conditions are given on $x$-axis, and  the odd
point-to-point reflection holds for (ii) and even for (v). The
author is not aware of any results for the conditions (iii) and
(iv). Our aim is to derive
 reflection formulas for the conditions (ii)-(v), that is, the formulas
expressing the value of a biharmonic function $u(x,y)$ at an
arbitrary point $(x_0 ,y_0)\in U_1$ in terms of its values at points
in $U_2$,
 when the
boundary conditions are given on a real-analytic curve, and to study
the properties of these mappings.


We remark that the structure of the reflection formulas attracts
attention of many researchers, interested, in particular, to answer  the
question: when the reflection is point to point \cite{khavinsons},
\cite{ek}, point to final set \cite{lopez} or point to continuous
set \cite{be-sa}.

The structure of the paper is as follows. In the section 2 we
discuss the so-called reflected fundamental solution for each case,
whose properties depend on the boundary condition and, therefore,
determine the structure of the reflection operator. In the section 3
we formulate and prove the main theorems.

\section{The reflected fundamental solutions}

This section is devoted to one of the key steps in deriving the
reflection formula, that is,  to the construction of the reflected fundamental
solutions for each case of the boundary conditions.

Let $\Gamma \subset \mathbb{R}^2$ be a non-singular real analytic curve
and point $P'\in \Gamma$. Consider the biharmonic differential operator in a neighborhood $U\subset \mathbb{R}^2$ of 
 the point $P'$. Let  $U_1$ and  $U_2$ be the components of $U\setminus \Gamma$, and 
 $P(x_0 ,y_0 )$ be a point in $U_1$ sufficiently close to $\Gamma$. The fundamental solution can be written in the form
$$
G=-\frac{1}{16\pi}((x-x_0 )^2 +(y-y_0 )^2 )\ln ((x-x_0 )^2 +(y-y_0
)^2 ) + g(x,y,x_0,y_0),
$$
where $g$ is a regular biharmonic function. It is obvious that $G$
is a real-analytic function in $\mathbb{R}^2$ except for the point $P(x_0 ,
y_0 )\in U_1$.

Let us complexify the problem, that is, consider a complex domain $W$ in the space $\mathbb{C}^2$ to which
the function $f$ defining the curve $\Gamma :=\{f(x,y)=0\}$  can be
continued analytically such that $W\cap \mathbb{R}^2=U$. Using the
change of variables $z=x+iy$, $w=x-iy,$ the equation of the
complexified curve $\Gamma _{\mathbb{C}}$ can be rewritten in the
form
$$
f\left ( \frac{z+w}{2}, \, \frac{z-w}{2i} \right ) =0,
$$
and if $grad \, f(x,y)\ne 0$ on $\Gamma$, can be also rewritten in
terms of the Schwarz function and its inverse, $w=S(z)$ and
$z=\widetilde{S}(w)$ \cite{davis}. In the characteristic variables
$G(z,w,z_0 ,w_0 )$ can be rewritten as
\begin{equation}\label{E1:1.4}
G  =  -\frac{1}{16\pi}( (z-z_0 )(w-w_0 )\ln [(z-z_0
)(w-w_0 )]+g(z,w,z_0,w_0).\nonumber
\end{equation}
It is obvious that the continuation of $G$  to the complex space (\ref{E1:1.4})
 has logarithmic singularities on the complex
characteristics passing through this point, i.e., on
$K_P:=\{(x-x_0)^2 + (y-y_0)^2=0\}.$

Note, that the specific choice of the regular part $g$ of the
fundamental solution $G$ does not affect our final result. Thus, for
 convenience  we choose the fundamental solution in the form
\begin{eqnarray}\label{E1:1.5}
& G  (z,w,z_0 ,w_0 )=  -\frac{1}{16\pi}(G_1 (z,w,z_0 ,w_0 ) +
G_2  (z,w,z_0 ,w_0 )) ,  \nonumber \\
& G_1  =  (z-z_0 )(w-w_0 )\bigl (\ln (z-z_0 ) -1\bigr ), \\
 &G_2
=(z-z_0 )(w-w_0 )\bigr ( \ln (w-w_0 )-1 \bigl ).\nonumber
\end{eqnarray}

Our goal is to construct  (multiple-valued) functions $\widetilde
G^{(j)}$, $j=\overline{1,5}$, which are biharmonic functions
satisfying on $\Gamma _{\mathbb{C}}$ one of the pair of conditions
(i) $\widetilde{G}=G$ and $\partial _n \widetilde{G}=\partial _n G$,
(ii) $\widetilde{G}=G$ and $\Delta \widetilde{G}=\Delta G$,
 (iii)
$\partial _n \widetilde{G}=\partial _n G$  and $\Delta
\widetilde{G}= \Delta G$, (iv) $\widetilde{G}= G$  and $\partial
_n\Delta \widetilde{G}=\partial _n \Delta G$, (v) $\partial _n
\widetilde{G}=\partial _n G$  and $\partial _n\Delta
\widetilde{G}=\partial _n \Delta G$, and having singularities only
on the characteristic lines intersecting the real space at point
$Q=R(P)$ (see formula (\ref{E:0.2})) in the domain $U_2$ and intersecting $\Gamma_{\mathbb{C}}$
at $K_P\cap \Gamma_{\mathbb{C}}$. These functions are called the
reflected fundamental solutions. According to (\ref{E1:1.5}) it is
convenient to seek these functions in the form
\begin{equation}\label{E:fund}
\widetilde{G}^{(j)}  (z,w,z_0 ,w_0 )=
-\frac{1}{16\pi}(\widetilde{G}_1^{(j)} (z,w,z_0 ,w_0 ) +
\widetilde{G}_2^{(j)} (z,w,z_0 ,w_0 )).
\end{equation}

It is easy to check that for the case (i) functions
$\widetilde{G}_1^{(1)}$ and $\widetilde{G}_2^{(1)}$ are
\begin{eqnarray}\label{E:1.4}
\widetilde{G}_1^{(1)}  = & (z-z_0 )(w-w_0 )\bigl (\ln
(\widetilde{S}(w)-z_0 ) -1\bigr
)+(z-\widetilde{S}(w))(w-w_0), \nonumber \\
 \widetilde{G}_2^{(1)} = & (z-z_0 )(w-w_0 )\bigr ( \ln (S(z)-w_0 )-1
\bigl )+ (w-S(z))(z-z_0).\nonumber
\end{eqnarray}
Similar functions (up to the regular part)  were used in \cite{as},
where function $g$ in the expression for $G$ (\ref{E1:1.4}) was
different from the chosen in (\ref{E1:1.5}).

To obtain $\widetilde{G}_1^{(j)}$ and $\widetilde{G}_2^{(j)}$ for
the cases
 (ii) -- (v),  we use the asymptotic expansions introduced by
Ludwig \cite{ludwig}, seeking these functions  in the form
\begin{equation}\label{Gja}
\widetilde{G}_1^{(j)}=\sum\limits_{k=1}^{\infty}b_k^{(j)}(z,w,w_0)f_k(\widetilde{S}(w)-z_0),
\end{equation}
\begin{equation}\label{Gjb}
\widetilde{G}_2^{(j)}=\sum\limits_{k=1}^{\infty}a_k^{(j)}(z,w,z_0)f_k
(S(z)-w_0),
\end{equation}
where
\begin{eqnarray}\label{fk}
 f_k(\xi )= \frac{\xi ^k}{k!}(\ln \xi -C_k),&\qquad k=0,1,...\, , \nonumber\\
 f_k(\xi )= (-1)^{-k-1}(-k-1)!\xi ^k,&\qquad k\le -1, C_0=0,\\
 C_k=  \sum\limits _{l=1}^k \frac{1}{l},&\qquad k=1,2,...\,
.\nonumber
\end{eqnarray}
This form ensures the location of singularities on the desired
characteristics. Substituting (\ref{Gja}) and (\ref{Gjb}) into the
biharmonic equation, we obtain the following differential equations
for the coefficients $a_k^{(j)}$ and $b_k^{(j)}$,
\begin{equation}
\frac{\partial ^2a_k^{(j)}}{\partial w^2}=0,\qquad\frac{\partial
^2b_k^{(j)}}{\partial z^2}=0,\qquad k=1,2,...\, .
\end{equation}

To find coefficients $a_k^{(j)}$ and $b_k^{(j)}$ for each case we
use the conditions on $\Gamma _{\mathbb{C}}$.

{\t Case (ii)}, $\widetilde{G}^{(2)}=G$ and $\Delta
\widetilde{G}^{(2)}=\Delta G$ on $\Gamma _{\mathbb{C}}$,
\begin{eqnarray}
a_1^{(2)}=z-z_0, &  \frac{\partial a_1^{(2)}}{\partial w}S^{\prime}=1
 \quad \mbox{ on } \Gamma _{\mathbb{C}}, \nonumber \\
a_k^{(2)}=0, & \qquad \qquad\qquad \frac{\partial
a_k^{(2)}}{\partial w}S^{\prime}=-\frac{\partial ^2
a_{k-1}^{(2)}}{\partial w\partial z}, \quad k=2,3,...\, ,\nonumber
\end{eqnarray}
\begin{eqnarray}
b_1^{(2)}=w-w_0, &  \frac{\partial b_1^{(2)}}{\partial z}\widetilde{S}^{\prime}=1
 \quad \mbox{ on } \Gamma _{\mathbb{C}}, \nonumber \\
b_k^{(2)}=0, & \qquad \qquad\qquad \frac{\partial
b_k^{(2)}}{\partial z}\widetilde{S}^{\prime}=-\frac{\partial ^2
b_{k-1}^{(2)}}{\partial w\partial z}, \quad k=2,3,...\, .\nonumber
\end{eqnarray}
Thus,
\begin{equation}\label{216}
a_1^{(2)}=z-z_0+\frac{(w-S(z))}{S^{\prime}}, \qquad
a_k^{(2)}=(-1)^{k-1}\frac{(w-S(z))}{S^{\prime}}\hat D^{k-1}_z,
\end{equation}
\begin{equation}\label{217}
b_1^{(2)}=w-w_0+\frac{(z-\widetilde{S}(w))}{\widetilde{S}^{\prime}},
\qquad
b_k^{(2)}=(-1)^{k-1}\frac{(z-\widetilde{S}(w))}{\widetilde{S}^{\prime}}\hat
D^{k-1}_w,
\end{equation}
where
$\hat D_z=\frac{\partial}{\partial z}\Bigl ( \frac{1}{S^{\prime}}  \Bigr )$,
$\hat D^2_z=\frac{\partial}{\partial z}\Bigl ( \frac{1}{S^{\prime}} \frac{\partial}{\partial z}\bigl ( \frac{1}{S^{\prime}}  \bigr ) \Bigr )$, etc, and
$\hat D_w=\frac{\partial}{\partial w}\Bigl ( \frac{1}{\widetilde{S}^{\prime}}  \Bigr )$.

Note that for the special case when $\Gamma$ is a straight line $y=0$, we have
\begin{equation}\label{Gja2}
\widetilde{G}_1^{(2)}=(z-w_0)(w-z_0)(\ln (w-z_0)-1 ),
\end{equation}
\begin{equation}\label{Gjb2}
\widetilde{G}_2^{(2)}=(w-z_0)(z-w_0)(\ln (z-w_0)-1).
\end{equation}

{\it Case (iii)}, $\partial _n \widetilde{G}^{(3)}=\partial _n G$
and $\Delta \widetilde{G}^{(3)}= \Delta G$ on $\Gamma
_{\mathbb{C}}$,
\begin{eqnarray}
a_1^{(3)}=-(z-z_0), &  \frac{\partial a_1^{(3)}}{\partial w}S^{\prime}=1
\quad \mbox{ on } \Gamma _{\mathbb{C}}, \nonumber \\
a_2^{(3)}S^{\prime}=1-\frac{\partial a_1^{(3)}}{\partial
z}+\frac{\partial a_1^{(3)}}{\partial w}S^{\prime}   , &
 \quad  \frac{\partial a_2^{(3)}}{\partial w}S^{\prime}=
 -\frac{\partial ^2 a_{1}^{(3)}}{\partial w\partial z}  , \nonumber \\
a_k ^{(3)} S^{\prime}=-\frac{\partial a_{k-1}^{(3)}}{\partial
z}+\frac{\partial a_{k-1}^{(3)}}{\partial w}S^{\prime}, & \quad
\qquad \frac{\partial a_k^{(3)}}{\partial w}S^{\prime}
=-\frac{\partial ^2 a_{k-1}^{(3)}}{\partial w\partial z} , \quad
k\ge 3 ,\nonumber
\end{eqnarray}
\begin{eqnarray}
b_1^{(3)}=-(w-w_0), &  \frac{\partial b_1^{(3)}}{\partial
z}\widetilde{S}^{\prime}=1
\quad \mbox{ on } \Gamma _{\mathbb{C}}, \nonumber \\
b_2^{(3)}\widetilde{S}^{\prime}=1-\frac{\partial b_1^{(3)}}{\partial
w}+ \frac{\partial b_1^{(3)}}{\partial z}\widetilde{S}^{\prime}  , &
\quad \frac{\partial b_2^{(3)}}{\partial z}\widetilde{S}^{\prime}=
-\frac{\partial ^2 b_{1}^{(3)}}{\partial w\partial z},\nonumber \\
b_k^{(3)}\widetilde{S}^{\prime}=-\frac{\partial
b_{k-1}^{(3)}}{\partial w}+\frac{\partial b_{k-1}^{(3)}}{\partial
z}\widetilde{S}^{\prime} , & \quad\qquad \frac{\partial
b_k^{(3)}}{\partial z}\widetilde{S}^{\prime} =-\frac{\partial ^2
b_{k-1}^{(3)}}{\partial w\partial z} , \quad k\ge 3 ,\nonumber
\end{eqnarray}
\begin{equation}
a_1^{(3)}=-(z-z_0)+\frac{(w-S(z))}{S^{\prime}}, \quad
b_1^{(3)}=-(w-w_0)+\frac{(z-\widetilde{S}(w))}{\widetilde{S}^{\prime}},
\end{equation}
\begin{equation}
a_k^{(3)}=(-1)^{k}\Bigl ( \frac{2k}{S^{\prime}}\hat D^{k-2}_z
-\frac{(w-S(z))}{S^{\prime}}\hat D^{k-1}_z \Bigr ), \quad k\ge 2,
\end{equation}
\begin{equation}
b_k^{(3)}=(-1)^{k}\Bigl (
\frac{2k}{\widetilde{S}^{\prime}}\hat
D^{k-2}_w-\frac{(z-\widetilde{S}(w))}{\widetilde{S}^{\prime}}\hat
D^{k-1}_w\Bigr ), \quad k\ge 2,
\end{equation}
where $\hat D^{0}_z=\hat D^{0}_w=1$.

For the special case when $\Gamma$ is a straight line $y=0$, we have
\begin{eqnarray}\label{Gja3}
\widetilde{G}_1^{(3)}=-(2w-z-w_0)(w-z_0)(\ln (w-z_0)-1 )\nonumber \\
+2(w-z_0)^2(\ln (w-z_0)-3/2 ),
\end{eqnarray}
\begin{eqnarray}\label{Gjb2}
\widetilde{G}_2^{(3)}=-(2z-w-z_0)(z-w_0)(\ln (z-w_0)-1)\nonumber \\
+2(z-w_0)^2(\ln (z-w_0)-3/2).
\end{eqnarray}

{\it Case (iv)}, $\widetilde{G}^{(4)}= G$  and $\partial _n\Delta
\widetilde{G}^{(4)}=\partial _n \Delta G$ on $\Gamma _{\mathbb{C}}$,
\begin{eqnarray}
a_1^{(4)}=z-z_0, &  \frac{\partial a_1^{(4)}}{\partial
w}S^{\prime}=-1
\quad \mbox{ on } \Gamma _{\mathbb{C}}, \nonumber \\
a_2^{(4)}=0, & \qquad  \frac{\partial a_2^{(4)}}{\partial w}(S^{\prime})^2=
-2\frac{\partial ^2 a_{1}^{(4)}}{\partial w\partial z}S^{\prime}
-\frac{\partial a_1^{(4)}}{\partial w}S^{\prime\prime} ,\nonumber \\
a_k^{(4)}=0,  \quad \frac{\partial a_k^{(4)}}{\partial
w}(S^{\prime})^2 & =-2\frac{\partial ^2 a_{k-1}^{(4)}}{\partial
w\partial z}S^{\prime} -\frac{\partial a_{k-1}^{(4)}}{\partial
w}S^{\prime\prime}- \frac{\partial ^3 a_{k-2}^{(4)}}{\partial
w\partial z^2} , \quad k\ge 3 ,\nonumber
\end{eqnarray}
\begin{eqnarray}
b_1^{(4)}=w-w_0, &  \frac{\partial b_1^{(4)}}{\partial
z}\widetilde{S}^{\prime}=-1
\quad \mbox{ on } \Gamma _{\mathbb{C}},\nonumber\\
b_2^{(4)}=0, & \qquad \frac{\partial b_2^{(4)}}{\partial z}(\widetilde{S}^{\prime})^2=
-2\frac{\partial ^2 b_{1}^{(4)}}{\partial w\partial z}\widetilde{S}^{\prime}
-\frac{\partial b_1^{(4)}}{\partial z}\widetilde{S}^{\prime\prime}  ,\nonumber \\
b_k^{(4)}=0,  \quad \frac{\partial b_k^{(4)}}{\partial
z}(\widetilde{S}^{\prime})^2 & =-2\frac{\partial ^2
b_{k-1}^{(4)}}{\partial w\partial z}\widetilde{S}^{\prime}
-\frac{\partial b_{k-1}^{(4)}}{\partial
z}\widetilde{S}^{\prime\prime}- \frac{\partial ^3
b_{k-2}^{(4)}}{\partial z\partial  w^2} , \quad k\ge 3 ,.\nonumber
\end{eqnarray}
\begin{equation}
a_1^{(4)}=z-z_0-\frac{(w-S(z))}{S^{\prime}},  \quad
a_k^{(4)}=(-1)^{k}\frac{(w-S(z))}{S^{\prime}}\hat D^{k-1}_z,
\end{equation}
\begin{equation}
b_1^{(4)}=w-w_0-\frac{(z-\widetilde{S}(w))}{\widetilde{S}^{\prime}},
\quad
b_k^{(4)}=(-1)^{k}\frac{(z-\widetilde{S}(w))}{\widetilde{S}^{\prime}}\hat
D^{k-1}_w, \,\, k\ge 2.
\end{equation}
For the special case when $\Gamma$ is a straight line $y=0$, we
obtain
\begin{equation}\label{Gja3}
\widetilde{G}_1^{(4)}=(2w-z-w_0)(w-z_0)(\ln (w-z_0)-1 ),
\end{equation}
\begin{equation}\label{Gjb2}
\widetilde{G}_2^{(4)}=(2z-w-z_0)(z-w_0)(\ln (z-w_0)-1).
\end{equation}

{\it Case (v)}, $\partial _n \widetilde{G}^{(5)}=\partial _n G$  and
$\partial _n\Delta \widetilde{G}^{(5)}=\partial _n \Delta G$ on
$\Gamma _{\mathbb{C}}$,
\begin{eqnarray}
a_1^{(5)}=-(z-z_0), &  \frac{\partial a_1^{(5)}}{\partial
w}S^{\prime}=-1
\quad \mbox{ on } \Gamma _{\mathbb{C}},\nonumber \\
a_2^{(5)}S^{\prime}=1-\frac{\partial a_1^{(5)}}{\partial z}+
\frac{\partial a_1^{(5)}}{\partial w}S^{\prime}, & \quad
\frac{\partial a_2^{(5)}}{\partial w}(S^{\prime})^2=
-2\frac{\partial ^2 a_{1}^{(5)}}{\partial w\partial z}S^{\prime}
-\frac{\partial a_1^{(5)}}{\partial w}S^{\prime\prime} ,\nonumber\\
a_k^{(5)}  S^{\prime}= -\frac{\partial a_{k-1}^{(5)}}{\partial z}
+\frac{\partial a_{k-1}^{(5)}}{\partial w}S^{\prime},& \nonumber \\
 \frac{\partial a_k^{(5)}}{\partial w}(S^{\prime})^2  =
 -2\frac{\partial ^2 a_{k-1}^{(5)}}{\partial w\partial z}S^{\prime} &
 -\frac{\partial a_{k-1}^{(5)}}{\partial w}S^{\prime\prime}-
\frac{\partial ^3 a_{k-2}^{(5)}}{\partial w\partial  z^2} , \quad
k\ge 3 ,\nonumber
\end{eqnarray}
\begin{eqnarray}
b_1^{(5)}=-(w-w_0), &  \frac{\partial b_1^{(5)}}{\partial
z}\widetilde{S}^{\prime}=-1
\quad \mbox{ on } \Gamma _{\mathbb{C}},\nonumber\\
b_2^{(5)}\widetilde{S}^{\prime}=1-\frac{\partial b_1^{(5)}}{\partial
w} +\frac{\partial b_1^{(5)}}{\partial z}\widetilde{S}^{\prime} , &
\qquad \frac{\partial b_2^{(5)}}{\partial
z}(\widetilde{S}^{\prime})^2= -2\frac{\partial ^2
b_{1}^{(5)}}{\partial w\partial z}\widetilde{S}^{\prime}
-\frac{\partial b_1^{(5)}}{\partial z}\widetilde{S}^{\prime\prime}  ,\nonumber \\
b_k^{(5)}\widetilde{S}^{\prime}=-\frac{\partial
b_{k-1}^{(5)}}{\partial w}
+\frac{\partial b_{k-1}^{(5)}}{\partial z}\widetilde{S}^{\prime}, & \nonumber \\
\frac{\partial b_k^{(5)}}{\partial z}(\widetilde{S}^{\prime})^2
=-2\frac{\partial ^2 b_{k-1}^{(5)}}{\partial w\partial
z}\widetilde{S}^{\prime} & -\frac{\partial b_{k-1}^{(5)}}{\partial
z}\widetilde{S}^{\prime\prime}- \frac{\partial ^3
b_{k-2}^{(5)}}{\partial z\partial  w^2} , \quad k\ge 3 ,.\nonumber
\end{eqnarray}
\begin{equation}
a_1^{(5)}=-(z-z_0)-\frac{(w-S(z))}{S^{\prime}},  \quad
\end{equation}
\begin{equation}
a_k^{(5)}=(-1)^{k-1}\Bigl (\frac{(2k-4)}{S^{\prime}}\hat
D^{k-2}_z-\frac{(w-S(z))}{S^{\prime}}\hat D^{k-1}_z\Bigr ), \quad
k\ge 2 ,
\end{equation}
\begin{equation}
b_1^{(5)}=-(w-w_0)-\frac{(z-\widetilde{S}(w))}{\widetilde{S}^{\prime}},
\end{equation}
\begin{equation}
b_k^{(5)}=(-1)^{k-1}\Bigl (\frac{(2k-4)}{\widetilde{S}^{\prime}}\hat
D^{k-2}_w- \frac{(z-\widetilde{S}(w))}{\widetilde{S}^{\prime}}\hat
D^{k-1}_w\Bigr ), \quad k\ge 2.
\end{equation}
For the special case when $\Gamma$ is a straight line $y=0$, this
imply
\begin{equation}\label{Gja3}
\widetilde{G}_1^{(5)}=-(z-w_0)(w-z_0)(\ln (w-z_0)-1 ),
\end{equation}
\begin{equation}\label{Gjb2}
\widetilde{G}_2^{(5)}=-(w-z_0)(z-w_0)(\ln (z-w_0)-1).
\end{equation}

Thus, we have constructed the reflected fundamental solutions for each
case of the boundary conditions (ii)-- (v) as formal series.
Convergence of the series follows from convergence of the
multipliers of the logarithms in (\ref{Gja}) and (\ref{Gjb}),
\begin{equation}\label{Vja}
{V}_1^{(j)}=\sum\limits_{k=1}^{\infty} b_k^{(j)}(z,w,w_0)
\frac{(\widetilde{S}(w)-z_0)^k}{k!},
\end{equation}
\begin{equation}\label{Vjb}
{V}_2^{(j)}=\sum\limits_{k=1}^{\infty} a_k^{(j)}(z,w,z_0)
\frac{(S(z)-w_0)^k}{k!}.
\end{equation}
Coefficients $a_k$ in (\ref{Vjb}) do not depend on $w_0$, therefore,
this expression can be interpreted as the Taylor series of
$V_2^{(j)}$ as a function of $-w_0$ at the point $-S(z)$. Function
$V_2^{(j)}$ can be described as the unique solution to the
Cauchy-Goursat problem for biharmonic functions with holomorphic
data on $\Gamma _{\mathbb C}$ and the characteristic line
$S(z)=w_0$. The data is prescribed by the boundary conditions for
each $j=\overline{2,5}$. For example, for $j=2$  the corresponding
problem is
\begin{eqnarray}
\Delta ^2 V_2^{(2)}  =0,& \mbox{ in } W\nonumber \\
V_2^{(2)}  =(z-z_0)(w-w_0)& \mbox{ on } \Gamma_{\mathbb C},\nonumber\\
\Delta V_2^{(2)}  =1& \mbox{ on } \Gamma_{\mathbb C},\\
V_2^{(2)} =0& \qquad\mbox{ on } S(z)-w_0=0 \nonumber,\\
\pd{V_2^{(2)}}{z} =(z-z_0)S^{\prime}+(w-S(z))& \qquad\mbox{ on }
S(z)-w_0=0.\nonumber
\end{eqnarray}
Note, that $S(z)$ is an analytic function in the neighborhood of
$\Gamma _{\mathbb C}$ and its derivative does not vanish on $\Gamma
_{\mathbb C}$. Existence and uniqueness of holomorphic solutions to
the Cauchy and Goursat problems for holomorphic partial differential
equations with holomorphic data are discussed in \cite{dima}.

The reflected fundamental solutions as well as functions
${V}_1^{(j)}$ and ${V}_2^{(j)}$ are used in the next section for
deriving the corresponding reflection formulas.

\section{The main result}

First we state the reflection formulas for biharmonic functions
given in the upper half plane describing the analytic continuation
across the the $x$-axis.

\noindent{\bf Theorem 3.1} {\it Let $U\subset \mathbb{R}^2$ be a sufficiently small domain divided by 
a straight line $\Gamma _0\subset U  :=\{y=0\}$ into two components $U_1\subset \mathbb{R}^2 _+$ and $U_2\subset \mathbb{R}^2 _-$,
$(x_0 ,y_0 )\in U_1$ and $(x_0,-y_0)\in U_2$. Then  any biharmonic function $u(x,y)$ in the domain $U_1$ 
subject to one of the
conditions (i)-(v) on  $\Gamma _0$ can be continued to the domain $U_2$, using the following formulas:

(i) If $u=\partial _n u =0$ on $\Gamma _0$, then
\begin{equation} \label{E:0.001p}
 u(x_0 ,y_0 ) =-u(x_0 ,-y_0 ) -2y_0 \pd{u}{y}(x_0 ,-y_0 )
-y_0 ^2 \Delta _{x,y} u(x_0 ,-y_0 ),
\end{equation}
(ii) if $u=\Delta u=0$ on $\Gamma _0$, then
\begin{equation} \label{E:0.002p}
 u(x_0 ,y_0 ) =-u(x_0 ,-y_0 ),
\end{equation}
(iii) if $u=\partial _n \Delta u =0$ on $\Gamma _0$, then
\begin{equation}\label{E:0.003p}
u(x_0 ,y_0 ) =-u(x_0 ,-y_0 )-y_0 \int\limits _{0}^{-y_{0}}\Delta
u(x_0,y)\, dy,
\end{equation}
where the integral is computed along the segment parallel to
$y$-axis,

 (iv) if $\partial _nu=\Delta u=0$ on $\Gamma _0$, then
\begin{equation}\label{E:0.004p}
u(x_0 ,y_0 ) =u(x_0 ,-y_0 ) - \int\limits _{0}^{-y_{0}}y\,\Delta
u(x_0,y)\,dy,
\end{equation}
(v) if $\partial _n u=\partial _n\Delta u=0$ on $\Gamma _0$, then
\begin{equation} \label{E:0.005p}
 u(x_0 ,y_0 ) =u(x_0 ,-y_0 ).
\end{equation}
}

\begin{remark}
If a biharmonic function $u(x,y)$ is also a harmonic function, then
(\ref{E:0.002p}) and (\ref{E:0.003p}) coincide with the odd
continuation (\ref{E:0.1}), while formulas (\ref{E:0.004p}) and
(\ref{E:0.005p}) with the even continuation (\ref{E:0.1b}) for
harmonic functions.
\end{remark}


\noindent{\bf Theorem 3.2} {\it Let $U\subset \mathbb{R}^2$ be a sufficiently small domain divided by 
a non-singular real analytic curve $\Gamma$ into two parts $U_1$ and $U_2$.
Let also $P(x_0 ,y_0 )$ be a point in $U_1$ having its reflected point $Q(R(x_0,y_0))$ in $U_2$ (see (\ref{E:0.2})). Then  any biharmonic function $u(x,y)$ in the domain $U$ 
subject to one of the
conditions (i)-(v) on the curve $\Gamma \subset U$ can be continued across $\Gamma$, using the following reflection formulas:

(i) if $u=\partial _n u =0$ on $\Gamma $, then
\begin{eqnarray}
 & u(P) = - u(Q)  -   \Bigl (x_0 -\frac{  S(x_0 +iy_0 ) +\widetilde{S}(x_0
-iy_0 )
          }{2} \Bigr  )\pd{u}{x}(Q) \nonumber \\
      &-  \Bigl ( y_0  + \frac{  S(x_0 +iy_0 ) -\widetilde{S}(x_0  -iy_0 )}{2i}
     \Bigr ) \pd{u}{y}(Q) \nonumber \\
     &- \frac{1}{4} \bigl ( x^2 _0 +y^2 _0  - S(x_0
      +iy_0 )(x_0 +iy_0 )  \label{E:1.2} \\
& - \widetilde{S}(x_0  -iy_0 ) (x_0  -iy_0 )+S(x_0 +iy_0
)\widetilde{S}(x_0 -iy_0 ) \bigr )
 \Delta _{x,y}u(Q) \nonumber,
\end{eqnarray}
(ii) if $u=\Delta u=0$ on $\Gamma $, then
\begin{equation} \label{E:0.002}
 u(P) =-u(Q)+\hat\mathbb{K}_2,
\end{equation}
(iii) if $u=\partial _n \Delta u =0$ on $\Gamma $, then
\begin{equation}\label{E:0.003}
u(x_0 ,y_0 ) =-u(x_0 ,-y_0 )+\hat\mathbb{K}_3,
\end{equation}
(iv) if $\partial _nu=\Delta u=0$ on $\Gamma $, then
\begin{equation}\label{E:0.004}
u(x_0 ,y_0 ) =u(x_0 ,-y_0 ) +\hat\mathbb{K}_4,
\end{equation}
(v) if $\partial _n u=\partial _n\Delta u=0$ on $\Gamma $, then
\begin{equation} \label{E:0.005}
 u(x_0 ,y_0 ) =u(x_0 ,-y_0 )+\hat\mathbb{K}_5,
\end{equation}
where
\begin{eqnarray}
& \hat\mathbb{K}_j =  \frac{1}{8i}\int\limits _{\Gamma}^{Q}\left (
V^{(j)}\pd{\Delta  u}{y}-\Delta  u \pd{V^{(j)}}{y}+ \Delta
V^{(j)}\pd{u}{y}-u\pd{\Delta V^{(j)}}{y}
  \right ) dx\label{E:1.2v} \\
&  - \left ( V^{(j)}\pd{\Delta  u}{x}-\Delta  u \pd{V^{(j)}}{x}+
\Delta V^{(j)}\pd{u}{x}-u\pd{\Delta  V^{(j)}}{x}
        \right )  dy, \nonumber
\end{eqnarray}
  the integral is computed along an arbitrary path joining the
curve $\Gamma$ with the reflected point $Q$,
$V^{(j)}=V^{(j)}_1-V^{(j)}_2$, $j=\overline{2,5}$.

}

It is obvious that the Theorem 3.1 is a special case of the Theorem
3.2., thus we will prove only the later for the boundary conditions
(ii)-(v) (for the case (i) see \cite{as}).

For simplicity, we assume
that $\Gamma$ is an algebraic curve. Under this assumption, the Schwarz
function and its inverse are analytic in the whole plane
$\mathbb{C}$ except for finitely many algebraic singularities.

{\it Proof of the Theorem 3.2.} The main step of the proof is already
done by constructing the reflected fundamental solution for each
case of boundary conditions (see section 2). The rest of the proof
is based on the contour deformation in  Green's formula \cite{G}
and is similar to \cite{as}.

The Green's formula, expressing the value of biharmonic function at
a point $P$ via the values of this function on a contour
$\gamma\subset U_1$ surrounding the point $P$, is
\begin{eqnarray}\label{green1}
    u(P)
= & \int\limits _{\gamma} \left  (
G\pd{\Delta  u}{y}-\Delta  u \pd{G}{y}+
\Delta  G\pd{u}{y}-u\pd{\Delta  G}{y}
  \right ) dx \label{E:1.20} \\
 & - \left (
G\pd{\Delta  u}{x}-\Delta  u \pd{G}{x}+
\Delta  G\pd{u}{x}-u\pd{\Delta  G}{x}
        \right )  dy,\nonumber
\end{eqnarray}
where  $G=G(x,y,x_0 ,y_0 )$ is an arbitrary fundamental solution of
the bi-Laplacian. The most suitable one for what follows is
(\ref{E1:1.5}).

Since the integrand in (\ref{E:1.20}) is a closed form, the value of
the integral does not change while we deform the contour $\gamma$
homotopically. Thus, the goal is to deform the contour $\gamma$ from
the domain $U_1$ to the domain $U_2$ by deforming it first  to the
complexified curve $\Gamma _{\mathbb{C}}$. This part of the
deformation is possible if the point $P$ lies so close to the curve
$\Gamma$ that there exists a connected domain $\Omega\subset \Gamma
_{\mathbb{C}}$ such that $\Omega$ contains both points of
intersections of the characteristic lines passing through the point
$P$ and  $\Omega$ can be univalently projected onto a plane domain
(for details, see \cite{sss}). Thus, we can replace the contour
$\gamma$ in (\ref{green1}) with the contour
$\gamma\prime\subset\Omega$, which is homotopic to $\gamma$ in
$\mathbb{C}^2\setminus
\{{(x-x_0)^2+(y-y_0)^2=0}\}=:\mathbb{C}^2\setminus K_P$.

 Note that due to homogeneous boundary conditions (ii) -- (v)
 a half of the terms (different for each case) in the integrand of
 (\ref{E:1.20}),
 while integrating along  $\gamma\prime$,  vanishes. For example,
  formula
(\ref{E:1.20}) in the case (ii) can be rewritten in the form
\begin{equation}\label{E:1.200}
    u(P)
=  \int\limits _{\gamma\prime} \left  (
G\pd{\Delta  u}{y}+\Delta  G \pd{u}{y}
  \right ) dx
  - \left (
G\pd{\Delta  u}{x}+\Delta  G \pd{u}{x}
        \right )  dy.
\end{equation}
To deform the contour $\gamma\prime$ from $\Gamma _{\mathbb{C}}$ to
the real domain $U_2$ we  replace the fundamental solution with the
corresponding reflected fundamental solution $\widetilde{G}^{(j)}$,
$j=\overline{1,5}$ (see section 2). Since functions
$\widetilde{G}^{(j)}$ have singularities only on the characteristic
lines intersecting the real space at point $Q=R(P)$ in the domain
$U_2$ and intersecting $\Gamma_{\mathbb{C}}$ at $K_P\cap
\Gamma_{\mathbb{C}}$, we are able to deform contour $\gamma\prime $
from the complexified curve $\Gamma _{\mathbb{C}}$ to the real
domain $U_2$ without changing the value of the integral \cite{sss}.
As a result, we obtain
\begin{eqnarray}
    u(P)
= & \int\limits _{\widetilde{\gamma}} \left  (
\widetilde{G}^{(j)}\pd{\Delta u}{y}-\Delta u
\pd{\widetilde{G}^{(j)}}{y}+ \Delta
\widetilde{G}^{(j)}\pd{u}{y}-u\pd{\Delta \widetilde{G}^{(j)}}{y}
  \right ) dx  \nonumber \\
 & - \left (
\widetilde{G}^{(j)}\pd{\Delta u}{x}-\Delta u
\pd{\widetilde{G}^{(j)}}{x}+ \Delta
\widetilde{G}^{(j)}\pd{u}{x}-u\pd{\Delta \widetilde{G}^{(j)}}{x}
        \right )  dy,\label{E:1.51}
\end{eqnarray}
where $\widetilde{\gamma}\subset U_2$ is a contour  surrounding the
point $Q$ and having endpoints on the curve $\Gamma$. Formula
(\ref{E:1.51}) in the characteristic variables has the form,
\begin{eqnarray}
    u(P)
= & 4 i\int\limits _{\widetilde{\gamma}} \left  (
\widetilde{G}^{(j)}\frac{\partial ^3 u}{\partial z^2 \partial w} +
\frac{\partial ^2 \widetilde{G}^{(j)}}{\partial z \partial w}
\pd{u}{z} -u\frac{\partial ^3 \widetilde{G}^{(j)}}{\partial z^2
\partial w} - \frac{\partial ^2 u}{\partial z \partial w}
\pd{\widetilde{G}^{(j)}}{z}
  \right ) dz \nonumber \\
 & - \left (
\widetilde{G}^{(j)}\frac{\partial ^3 u}{\partial z \partial w^2 } +
\frac{\partial ^2 \widetilde{G}^{(j)}}{\partial z \partial w}
\pd{u}{w} -u\frac{\partial ^3 \widetilde{G}^{(j)}}{\partial z
\partial w^2 } - \frac{\partial ^2 u}{\partial z \partial w}
\pd{\widetilde{G}^{(j)}}{w}
        \right )  dw.\label{E:1.6}
\end{eqnarray}
If we substitute    (\ref{E:fund}), (\ref{Gja}) and (\ref{Gjb}) into
(\ref{E:1.6}) and move one endpoint of the contour
$\widetilde{\gamma}$ along the curve $\Gamma$ to the other endpoint,
the integral terms containing products of the function $u$ and
regular part of the function $\widetilde{G}^{(j)}$ and their
derivatives vanish, while the integral terms containing logarithms
can be combined as follows,

\begin{eqnarray}\label{E:1.066}
 & -\frac{i}{4\pi}\int\limits _{\widetilde{\gamma}}
 \ln (\widetilde{S}(w)-z_0 ) \Bigl \{
\bigl ( V_1^{(j)}\frac{\partial ^3 u}{\partial z^2 \partial w}+
  \frac{\partial ^2V_1^{(j)} }{\partial z\partial w} \pd{u}{z}
 -\frac{\partial ^2 u}{\partial z\partial w}
 \frac{\partial  V_1^{(j)}}{\partial z}\bigr ) dz \nonumber \\
& - \bigl ( V_1^{(j)}\frac{\partial ^3 u}{\partial z\partial w^2}+
\frac{\partial ^2V_1^{(j)} }{\partial z\partial w}\pd{u}{w}
-u\frac{\partial ^3 V_1^{(j)}}{\partial z \partial w^2}
-\frac{\partial ^2 u}{\partial z\partial w}\pd{V_1^{(j)}}{w}\bigr )
dw \Bigr \},
\end{eqnarray}
\begin{eqnarray}
 &-\frac{i}{4\pi}\int\limits _{\widetilde{\gamma}}\ln (S(z)-w_0 )
\Bigl \{ \bigl ( V_2^{(j)}\frac{\partial ^3 u}{\partial z^2 \partial
w}+
  \frac{\partial ^2V_2^{(j)} }{\partial z\partial w} \pd{u}{z}
  -u\frac{\partial ^3 V_2^{(j)}}{\partial z^2 \partial w}
 -\frac{\partial ^2 u}{\partial z\partial w}
 \frac{\partial  V_2^{(j)}}{\partial z}\bigr ) dz \nonumber \\
& - \bigl ( V_2^{(j)}\frac{\partial ^3 u}{\partial z\partial w^2}+
\frac{\partial ^2V_2^{(j)} }{\partial z\partial w}\pd{u}{w}
-\frac{\partial ^2 u}{\partial z\partial w}\pd{V_2^{(j)}}{w}\bigr )
dw \Bigr \} ,
\end{eqnarray}
where $\widetilde{\gamma} $ is a loop surrounding the point $Q$ and
having endpoints on the curve $\Gamma$.

The  logarithm  $\ln (\widetilde{S}(w)-z_0 )$ obtains the increment
$2\pi i$ along the loop, while the logarithm $\ln (S(z)-w_0 )$ gets $(-2\pi i)$. Thus, compressing $
\widetilde{\gamma}$ to a segment joining $Q$ to $\Gamma $, we find that the
integral  (\ref{E:1.066}) can be rewritten as
\begin{eqnarray}\label{E:1.066a}
\, & \hat\mathbb{K}_j= \frac{1}{2}\int\limits _{\Gamma}^{Q} \Bigl \{
\bigl ( V^{(j)}\frac{\partial ^3 u}{\partial z^2 \partial w}+
  \frac{\partial ^2V^{(j)} }{\partial z\partial w} \pd{u}{z}
  +u\frac{\partial ^3 V_2^{(j)}}{\partial z^2 \partial w}
 -\frac{\partial ^2 u}{\partial z\partial w}
 \frac{\partial  V^{(j)}}{\partial z}\bigr ) dz \nonumber \\
& - \bigl ( V^{(j)}\frac{\partial ^3 u}{\partial z\partial w^2}+
\frac{\partial ^2V^{(j)} }{\partial z\partial w}\pd{u}{w}
-u\frac{\partial ^3 V_1^{(j)}}{\partial z \partial w^2}
-\frac{\partial ^2 u}{\partial z\partial w}\pd{V^{(j)}}{w}\bigr ) dw
\Bigr \},
\end{eqnarray}
where $V^{(j)}=V_1^{(j)}-V_2^{(j)}$.

The rest of nonzero terms in (\ref{E:1.6}) are terms involving
derivatives of the logarithms. Some of these integrals also vanish
due to the properties of functions $V_1^{(j)}$ and $V_2^{(j)}$,
resulting in
\begin{eqnarray}\label{E:1.66c}
  &- \frac{i}{4\pi}   \int\limits _{\widetilde{\gamma}} \Bigl  (
- \pd{a_1^{(j)}}{w}\frac{ (S^{\prime}(z))^2 u}{ S(z) -w_0 }dz +
 \pd{b_1^{(j)}}{z}\frac{ ((\widetilde{S}^{\prime}(w))^2 u }{ \widetilde{S}(w) -z_0  }
        dw\Bigr )
  \nonumber \\
 & = -\frac{1}{2}u(Q)(\pd{a_1^{(j)}}{w}(Q)S^{\prime}(Q)
 +\pd{b_1^{(j)}}{z}(Q)\widetilde{S}^{\prime}(Q)).
\end{eqnarray}
Combining (\ref{E:1.066a}) and (\ref{E:1.66c})  we finally obtain,
\begin{equation}\label{E:1.25}
    u(P)=-\frac{1}{2}u(Q)(\pd{a_1^{(j)}}{w}(Q)S^{\prime}(Q)
 +\pd{b_1^{(j)}}{z}(Q)\widetilde{S}^{\prime}(Q))+\hat\mathbb{K}_j.
\end{equation}
Here the expression in the parentheses with the appropriate choice
of the coefficients $a_1^{(j)}$ and $b_1^{(j)}$ (see section 2) is
equal to either 1 or -1. Thus, formula (\ref{E:1.25}) in variables $x,
y$ is equivalent to (\ref{E:0.002}) -- (\ref{E:1.2v}).

\begin{remark}
Formula (\ref{E:1.25}) gives continuation of a biharmonic function
from the domain $U_1\subset\mathbb{R}^2$ to the domain
$U_2\subset\mathbb{R}^2$ as a multiple-valued function whose
singularities coincide with the singularities of the functions $S$
or $\widetilde{S}$.
\end{remark}

\begin{remark}
In the special case when $\Gamma$ is a straight line and boundary
conditions (ii) or (v) are applied, $V_1^{(j)}=V_2^{(j)}$, and,
therefore, $\hat\mathbb{K}_j\equiv 0$, $j=2,5$.
\end{remark}

\begin{example}
As an example of formula (\ref{E:1.25}) consider a biharmonic function $u(x,y)$ subject to the Navier conditions (ii), $u=\Delta u=0$, on a unit circle centered at the origin, $x^2+y^2=1$.
In this case series (\ref{Vja}) and (\ref{Vjb}) with coefficients (\ref{216}), (\ref{217}) can be summed, and $V^{(2)}_2=(z-z_0)(1/z-w_0)+(1/z-w)(1/w_0-z)$. The reflection formula then has the form 
\begin{equation}\label{circle}
u(r_0,\theta_0)=-u(\frac{1}{r_0},\theta_0) +
\frac{r_0^2-1}{4r_0}\int\limits_{1}^{r_0^{-1}}
\frac{1-r^2}{r^2}\, \Bigl ( \frac{1}{r} \partial _r u(r,\theta _0)+\partial ^2 _{rr}u(r,\theta _0)\Bigr )\, dr,
\end{equation}
where the integral is computed along the straight line $\theta
=\theta _0$.
\end{example}

\textbf{Acknowledgments.} I would like to thank Professor Reinhard Farwig for
 sending me reprints of his interesting papers.

\providecommand{\bysame}{\leavevmode\hbox to3em{\hrulefill}\thinspace}

\end{document}